
\input amstex
\documentstyle{amsppt}
\magnification =\magstep 1
\nologo
\NoBlackBoxes
\baselineskip =18pt
\pagewidth{5.5 in}
\vsize =9truein

\topmatter

\title
Discrete Cocompact subgroups of $\roman G_{5,3}$ \\\ 
and related $C^{*}$-algebras
\endtitle

\rightheadtext{Discrete cocompact subgroups of $\roman G_{5,3}$}

\author P.~Milnes and S.~Walters \endauthor

\address
Department of Mathematics, University of Western Ontario, London,
Ontario N6A 5B7, Canada
\endaddress
\email milnes\@uwo.ca \hfill 
\endemail

\address
Department of Mathematics and Computer Science, University of
Northern British Columbia, Prince George, British Columbia V2N 4Z9, Canada
\endaddress
\email walters\@hilbert.unbc.ca  \hfill 
\newline 
\indent {\it Home page}: http://hilbert.unbc.ca/walters
\endemail

\thanks Research partly supported by NSERC
grants A7857 and OGP0169928 \endthanks

\date
\sevenrm February 17, 2001
\enddate



\abstract
The discrete cocompact subgroups of the 5-dimensional Lie group $\roman G_{
5,3}$ are determined
up to isomorphism.  Each of their group $C^{*}$-algebras is studied by determining all of 
its simple infinite dimensional quotient $C^{*}$-algebras.
The $K$-groups and trace invariants of the latter are also obtained.
\endabstract

\endtopmatter


\def\vartau{\tau\hskip2pt\!\!\!\!\iota\hskip1pt}


\def\EH{1}   
\def\HF{2}   
\def\AM{3}
\def\MWa{4}
\def\MWb{5}
\def\MWc{6}
\def\ON{7}
\def\OP{8}
\def\JPa{9}
\def\JPb{10}
\def\GP{11}
\def\MP{12}
\def\PV{13}
\def\MR{14}
\def\SW{15}
\def\ZM{16}  


\document

\subhead \S 1. Introduction \endsubhead

Consider the Lie group $\roman G_{5,3}$ equal to $\Bbb R^5$ as a set with
multiplication given by
$$(h,j,k,m,n)(h^{\prime},j^{\prime},k^{\prime},m^{\prime},n^{\prime}
)=$$
$$(h+h^{\prime}+nj^{\prime}+m^{\prime}n(n-1)/2+mk^{\prime},j+j^{\prime}
+nm^{\prime},k+k^{\prime},m+m^{\prime},n+n^{\prime}).$$
and inverse 
$$(h,j,k,m,n)^{-1}=(-h+nj+mk-mn(n-1)/2,-j+nm,-k,-m,-n).$$
The group $\roman G_{5,3}$ is one of only six
nilpotent, connected, simply connected, 5-dimensional Lie groups; it
seemed the most tractable of them for our present purposes. (Our
notation is as in
Nielsen [\ON ], where a detailed catalogue of Lie groups like this one 
is given.) In [\MWb , Section 3] the authors
have studied a natural discrete cocompact subgroup $\roman H_{5,3}$, the lattice 
subgroup $\roman H_{5,3}=\Bbb Z^5\subset \roman G_{5,3}$.  In section 2 of this paper we 
study the group $\roman G_{5,3}$ more closely, determining the 
isomorphism classes of all its discrete cocompact subgroups 
(Theorem 1).  These are given
by five integer parameters $\alpha ,\beta ,\gamma ,\delta ,\epsilon$ that satisfy certain
conditions (see $(*)$ and $(**)$ of Theorem 1), and are denoted by
$\roman H_{5,3}(\alpha ,\beta ,\gamma ,\delta ,\epsilon )$.
It is shown that each such subgroup is isomorphic to a cofinite subgroup of 
$\roman H_{5,3}=\roman H_{5,3}(1,0,1,1,0)$.
Conversely, each cofinite subgroup of $\roman H_{5,3}\subset \roman G_{
5,3}$ is a discrete 
cocompact subgroup of $\roman G_{5,3}$.  In sections 3 and 4 the group
$C^{*}$-algebras of the $\roman H_{5,3}(\alpha ,\beta ,\gamma ,\delta 
,\epsilon )$'s
are examined by obtaining their simple infinite dimensional quotients.
These are shown to be crossed products of certain types of Heisenberg
$C^{*}$-algebras (in Packer's terminology [\JPb ]) and the rest are matrix algebras over
irrational rotation algebras (Theorem 5).  In section 5 the $K$-groups of the simple
quotients are calculated (Theorem 6) as are their trace invariants (Theorem 8).
The paper ends with a discussion of the classification of the simple 
quotients.


\medpagebreak

We use the conventional notation for crossed products as in, for example, [\GP ] or [\ZM ].
Hence, if a discrete group $G$ acts on a $C^{*}$-algebra $A$,
we write $C^{*}(A,G)$ to denote the
associated $C^{*}$-crossed product algebra.  We use a similar notation for twisted crossed 
products, i.e. when there is a cocycle instead of an action (as in Theorem 2).
(See the Preliminaries of [\MWb ] for more details.)

\bigpagebreak

\subhead \S 2. Determination of the Discrete Cocompact Subgroups \endsubhead

\proclaim{1. Theorem} Every discrete cocompact subgroup $\roman H$ of $
\roman G_{5,3}$ has 
the following form$:$ there are integers $\alpha$, $\beta$, $\gamma$, $
\delta$ and 
$\epsilon$ satisfying $\alpha ,\gamma ,\delta >0$, and
$$0\leq\epsilon\leq \roman g\roman c\roman d\,\{\gamma ,\delta \}/
2\text{ and}\leqno (*)$$
$$0\leq\beta\leq \roman g\roman c\roman d\,\{\alpha ,\gamma ,\delta 
,\epsilon \}/2,\leqno (**)$$
yielding $\roman H\cong \roman H_{5,3}(\alpha ,\beta ,\gamma ,\delta 
,\epsilon )$ $(\,=\Bbb Z^5$ as a set$)$ with multiplication 
$$\left\{\aligned
&(h,j,k,m,n)(h^{\prime},j^{\prime},k^{\prime},m^{\prime},n^{\prime}
)=\\
&(h+h^{\prime}+\gamma nj^{\prime}+\alpha\gamma m^{\prime}n(n-1)/2+
\beta nm^{\prime}+\delta mk^{\prime}+\epsilon nk^{\prime},\\
&j+j^{\prime}+\alpha nm^{\prime},k+k^{\prime},m+m^{\prime},n+n^{\prime}
).\endaligned
\right.\leqno (\roman m)$$
Different choices for $\alpha$, $\beta$, $\gamma$, $\delta$ and $\epsilon$ give 
non-isomorphic groups.  Each such group is, in fact, isomorphic to a 
cofinite subgroup of 
$\roman H_{5,3}$ $($the lattice subgroup of $\roman G_{5,3})$, and each cofinite subgroup
of $\roman H_{5,3}$ is isomorphic to some $\roman H_{5,3}(\alpha ,
\beta ,\gamma ,\delta ,\epsilon )$.  
\endproclaim

\flushpar PROOF. Using the discreteness and cocompactness as in [\MWc ], 
the second commutator subgroup of H tells us that there is a member 
(with entries that don't need to be identified indicated by $*)$

$e_5=(*,*,*,\frak a,z)$ 

\flushpar of H, where $z>0$ is the smallest 
positive number that can appear as the last coordinate of a member 
of H. Continuing in this vein, we get 

$e_4=(*,*,*,y,0)$, 

$e_3=(*,\frak b,x,0,0)$, 

$e_2=(*,w,0,0,0)$ and 

$e_1=(v,0,0,0,0)$,

\flushpar where $x>0$ is the smallest 
positive number that can appear as the 3rd coordinate of a member of 
H whose last 2 coordinates are 0, and similarly for $v,w$ and $y$. Also, 
all other coordinates are $\geq 0$, and the bottom non-zero coordinate in 
each column is greater than the coordinates above it, e.g., $w>\frak b
\geq 0$
and $w$ is also greater than the 2nd coordinate of $e_5$ or of $e_
4$.
These considerations show that the map 
$$\pi :(h,j,k,m,n)\mapsto e_1^he_2^je_3^ke_4^me_5^n,\,\,\,\Bbb Z^5
\to \roman H,$$
is 1 -- 1 and onto. We want the multiplication (m) for $\Bbb Z^5$ that makes $
\pi$ a 
homomorphism (hence an isomorphism); (m) is determined using the 
commutators,
$$\left\{\aligned
&[e_5,e_4]=(*,zy,0,0,0)=e_1^{\beta}e_2^{\alpha},\,\,\,[e_5,e_3]=(z
\frak b+x\frak a,0,0,0,0)=e_1^{\epsilon},\\
&[e_5,e_2]=(zw,0,0,0)=e_1^{\gamma},\text{ and }[e_4,e_3]=(xy,0,0,0
)=e_1^{\delta},\endaligned
\right.\leqno (C)$$
for some integers $\alpha ,\beta ,\gamma ,\delta ,\epsilon$ (other pairs of $
e$'s commuting). Using the 
commutators to collect terms in
$$(e_1^he_2^je_3^ke_4^me_5^n)(e_1^{h^{\prime}}e_2^{j^{\prime}}e_3^{
k^{\prime}}e_4^{m^{\prime}}e_5^{n^{\prime}})$$
gives the multiplication formula (m) for $\Bbb Z^5$, and also the equation
$$e^n_5e^{m^{\prime}}_4=e_1^{\alpha\gamma m^{\prime}n(n-1)/2+\beta 
nm^{\prime}}e_2^{\alpha m^{\prime}n}e_4^{m^{\prime}}e_5^n,$$
which the reader may find helpful in checking computations later.

\bigpagebreak

For a start in putting the restrictions on $\alpha ,\beta ,\gamma 
,\delta ,\epsilon$, $(C)$ tells 
us that $\alpha ,\gamma ,\delta >0$ (since $v$, $w$, $x$, $y$ and $
z>0$). Let $Z$ denote the center
of $\roman H_{5,3}(\alpha ,\beta ,\gamma ,\delta ,\epsilon )$, $Z=
(\Bbb Z,0,0,0,0)$.
Then, as for $\roman G_4$, with quotients and subgroups it is shown that 
different (positive) $\alpha$, $\gamma$, $\delta$
give non-isomorphic groups, e.g., $\roman H_{5,3}(\alpha ,\beta ,\gamma 
,\delta ,\epsilon )/Z$ gives $\alpha$, then $Z$ 
modulo the second commutator subgroup gives $\gamma$, and also, 
with $K_3$, $K_4\subset \roman H_{5,3}(\alpha ,\beta ,\gamma ,\delta 
,\epsilon )$ as defined below,
$$Z\supset (\delta \Bbb Z,0,0,0,0)=\{xyx^{-1}y^{-1}\,\,\,|\,\,\,x\in 
K_3,\,\,\,y\in K_4\}$$
and $Z/(\delta \Bbb Z,0,0,0,0)=\Bbb Z_{\delta}$, the cyclic group of order $
\delta$.

Then we have an isomorphism 
of $\roman H_{5,3}(\alpha ,\beta ,\gamma ,\delta ,\epsilon )\text{ onto }
\roman H_{5,3}(\alpha ,\beta ,\gamma ,\delta ,\epsilon +d\gamma +e
\delta )$,
which is simpler to give in terms 
of generators,
$$e_3\mapsto e_3^{\prime}=e_2^de_3,\,\,\,e_5\mapsto e_5^{\prime}=e_
4^ee_5,\text{ and }e_i\mapsto e_i^{\prime}=e_i\text{ otherwise.}\leqno 
(\circledast )$$
Here we are merely changing the basis for $\roman H_{5,3}(\alpha ,
\beta ,\gamma ,\delta ,\epsilon )$, and the only 
commutator (using (m) and $(C)$) that changes is $[e_5^{\prime},e_
3^{\prime}]=e_1^{\epsilon +e\delta +d\gamma}$, so the 
resulting isomorphism is of $\roman H_{5,3}(\alpha ,\beta ,\gamma 
,\delta ,\epsilon )\text{ onto }\roman H_{5,3}(\alpha ,\beta ,\gamma 
,\delta ,\epsilon +d\gamma +e\delta ),$
which shows we can require
$$0\leq\epsilon <\roman g\roman c\roman d\,\{\gamma ,\delta \}.$$
This, accompanied by another isomorphism,
$$(h,j,k,m,n)\mapsto (-h,-j,k,-m,n),\,\,\,\roman H_{5,3}(\alpha ,\beta 
,\gamma ,\delta ,\epsilon )\to \roman H_{5,3}(\alpha ,\beta ,\gamma 
,\delta ,-\epsilon ),\leqno (\circledast^{\prime})$$
assures that we can have 
$$0\leq\epsilon\leq \roman g\roman c\roman d\,\{\gamma ,\delta \}/
2,\leqno (*)$$
the required range for $\epsilon$. 

Now, to control $\beta$,
$$\left\{\aligned
&e_1\mapsto e_1=e_1^{\prime},\,\,\,e_2\mapsto e_1^{-q}e_2=e_2^{\prime}
,\,\,\,e_3\to e_3=e_3^{\prime},\\
&e_4\mapsto e_2^re_3^ge_4\text{ and }e_5\mapsto e_3^{-f}e_5=e_5^{\prime}\endaligned
\right.\leqno (\dagger )$$
is an isomorphism of $\roman H_{5,3}(\alpha ,\beta ,\gamma ,\delta 
,\epsilon )$ onto $\roman H_{5,3}(\alpha ,\beta +q\alpha +r\gamma 
+f\delta +g\epsilon ,\gamma ,\delta ,\epsilon )$, which 
yields
$$0\leq\beta <\roman g\roman c\roman d\{\alpha ,\gamma ,\delta ,\epsilon 
\}.$$
Then the isomorphism 
$$(h,j,k,m,n)\mapsto (-h,j,k,-m,-n)\leqno (\dagger^{\prime})$$
of $\roman H_{5,3}(\alpha ,\beta ,\gamma ,\delta ,\epsilon )$ onto $
\roman H_{5,3}(\alpha ,-\beta +\alpha\gamma ,\gamma ,\delta ,\epsilon 
)$ leads to the conclusion 
$$0\leq\beta\leq \roman g\roman c\roman d\,\{\alpha ,\gamma ,\delta 
,\epsilon \}/2.\leqno (**)$$

It must still be shown that changing 
$\epsilon$ or $\beta$ within the allowed limits (namely, $\epsilon$ and $
\beta$ must satisfy
$(*)$ and $(**)$, respectively) gives a non-isomorphic group.

\smallpagebreak

So, suppose that $\varphi :\roman H=\roman H_{5,3}(\alpha ,\beta ,
\gamma ,\delta ,\epsilon )\to \roman H_{5,3}(\alpha ,\beta^{\prime}
,\gamma ,\delta ,\epsilon^{\prime})=\roman H^{\prime}$
is an isomorphism. Then 
$$\varphi :Z=K_1=(\Bbb Z,0,0,0,0)\to (\Bbb Z,0,0,0,0)=K_1^{\prime}
=Z^{\prime},$$
$$K_2=(\Bbb Z,\Bbb Z,0,0,0)\to (\Bbb Z,\Bbb Z,0,0,0)=K_2^{\prime},$$
$$K_3=(\Bbb Z,\Bbb Z,\Bbb Z,0,0)\to (\Bbb Z,\Bbb Z,\Bbb Z,0,0)=K_3^{
\prime},\text{ and}$$
$$K_4=(\Bbb Z,\Bbb Z,\Bbb Z,\Bbb Z,0)\to (\Bbb Z,\Bbb Z,\Bbb Z,\Bbb Z
,0)=K_4^{\prime},$$
since the $Z$'s are the centers, the $K_2$'s consist of those $s\in 
\roman H$ for 
which $s^r$ is in the commutator subgroup of H for some $r\in \Bbb Z$, the $
K_3$'s 
are the largest subsets for which all 
commutators are central (e.g., $xyx^{-1}y^{-1}\in Z$ for all $x\in 
K_3$ and 
$y\in \roman H_{5,3}(\alpha ,\beta ,\gamma ,\delta ,\epsilon )$),
and the $K_4$'s are 
the centralizers of the commutator subgroups. So
we must have 
$$\varphi (0,0,0,0,1)=(*,*,-f,e,a)=S_5\text{ with }a=\pm 1,$$
$$\varphi (0,0,0,1,0)=(*,r,g,b,0)=S_4\text{ with }b=\pm 1,\text{ and}$$
$$\varphi (0,0,1,0,0)=(*,d,c,0,0)=S_3\text{ with }c=\pm 1;$$
furthermore, commutators give
$$\varphi (\beta ,\alpha ,0,0,0)=[S_5,S_4]=S_5S_4S_5^{-1}S_4^{-1}=
(*,\alpha ab,0,0,0),$$
hence $\varphi (0,1,0,0,0)=(q,ab,0,0,0)=S_2$, and
$$\varphi (\gamma ,0,0,0,0)=[S_5,S_2]=(\gamma a^2b,0,0,0,0),$$
so $\varphi (1,0,0,0,0)=(b,0,0,0,0)=S_1$, but also
$$\varphi (\delta ,0,0,0,0)=[S_4,S_3]=(\delta bc,0,0,0,0),$$
so $c=1$. 
Furthermore, $\varphi (\epsilon ,0,0,0,0)=[S_5,S_3]=(a\epsilon^{\prime}
+e\delta +ad\gamma ,0,0,0,0)$, which 
shows that the manipulations at $(\circledast )$ and $(\circledast^{
\prime})$ above 
give the only way of changing $\epsilon$ in $\roman H_{5,3}(\alpha 
,\beta ,\gamma ,\delta ,\epsilon )$; that is, if 
$$0\leq\epsilon ,\epsilon^{\prime}\leq \roman g\roman c\roman d\,\{
\gamma ,\delta \}/2\leqno (*)$$
and $\epsilon =\pm\epsilon^{\prime}+a_1\delta +a_2\gamma$ with $a_
1,a_2\in \Bbb Z$, then $\epsilon =\epsilon^{\prime}$. Now consider
$$\varphi (h,j,k,m,n)=\varphi ((h,0,0,0,0)(0,j,0,0,0)(0,0,k,0,0)(0
,0,0,m,0)(0,0,0,0,n))$$
$$=(hS_1)\cdot (jS_2)\cdot (kS_3)\cdot S_4^m\cdot S_5^n=hS_1+jS_2+
kS_3+S_4^m\cdot S_5^n\in \roman H^{\prime}.$$
Note that $S_5^n\neq nS_5$, but $S_5^n=(*,*,-nf,ne,na)$, and also 
$S_4^m=(*,mr,mg,mb,0)$; further, the $(jS_2)$ term puts a $jq$ in the first 
entry of $\varphi (h,j,k,m,n)$, so also $(j+j^{\prime}+\alpha nm^{
\prime})q$ in the first entry of
$\varphi (h,j,k,m,n)\cdot\varphi (h^{\prime},j^{\prime},k^{\prime}
,m^{\prime},n^{\prime})$ (product in $\roman H_{5,3}(\alpha ,\beta^{
\prime},\gamma ,\delta ,\epsilon )$).
Then, equating the coefficients of the $nm^{\prime}$ terms in the first entry of 
$$\varphi (e_5^ne_4^{m^{\prime}})\text{ and }\varphi (e_5^n)\varphi 
(e_4^{m^{\prime}})=S_5^nS_4^{m^{\prime}}\text{ gives}$$
$$b(-\alpha\gamma /2+\beta )+q\alpha =ab\beta^{\prime}-ab\alpha\gamma 
/2+ag\epsilon +ar\gamma +(eg+bf)\delta ,\text{ or}$$
$$\beta =\pm\beta^{\prime}+a_1\alpha +a_2\gamma +a_3\delta +a_4\epsilon\text{ for some }
a_i\in \Bbb Z,\,\,\,1\leq i\leq 4,$$
which shows that the manipulations at $(\dagger )$ and $(\dagger^{
\prime})$ above give the only 
way of changing just $\beta$ in $\roman H_{5,3}(\alpha ,\beta ,\gamma 
,\delta ,\epsilon )$.

\smallpagebreak

Here is an isomorphism 
$\varphi$ of $\roman H_{5,3}(\alpha ,\beta ,\gamma ,\delta ,\epsilon 
)$ into the lattice subgroup 
$\roman H_{5,3}=\Bbb Z^5\subset \roman G_{5,3}$ in terms of generators; $
\roman H_{5,3}$ 
has multiplication 
$$\left\{\aligned
&(h,j,k,m,n)(h^{\prime},j^{\prime},k^{\prime},m^{\prime},n^{\prime}
)=\\
&(h+h^{\prime}+nj^{\prime}+m^{\prime}n(n-1)/2+mk^{\prime},\\
&j+j^{\prime}+nm^{\prime},k+k^{\prime},m+m^{\prime},n+n^{\prime})\endaligned
\right.\leqno (\roman m^{\prime})$$
(i.e., $\alpha =\gamma =\delta =1$ and $\beta =\epsilon =0$). First suppose $
\epsilon >0$.
Then, with $\frak d=\alpha\gamma\epsilon$ and generators
$$e_1=(1,0,0,0,0),\,\,\,e_2=(0,1,0,0,0),\,\,\,...\,\,\,,\,\,\,e_5=
(0,0,0,0,1)$$
for $\roman H_{5,3}(\alpha ,\beta ,\gamma ,\delta ,\epsilon )$ satisfying 
$$[e_5,e_4]=e_1^{\beta}e_2^{\alpha},\,\,\,[e_5,e_3]=e_1^{\epsilon}
,\,\,\,[e_5,e_2]=e_1^{\gamma},\text{ and }[e_4,e_3]=e_1^{\delta},\leqno 
(C)$$
$\varphi$ is given by 
$$\varphi :e_1\mapsto e_1^{\prime}=(\delta \frak d^2,0,0,0,0),\,\,\,
e_2\mapsto e_2^{\prime}=(\gamma\delta \frak d(\frak d-1)/2,\gamma\delta 
\frak d,0,0,0),$$
$$e_3\mapsto e_3^{\prime}=(0,\delta\epsilon \frak d,\delta\epsilon 
\frak d,0,0),\,\,\,e_4\mapsto e_4^{\prime}=(0,\beta\delta \frak d,
0,\alpha\gamma\delta ,0),$$
$$\text{ and }e_5\mapsto e_5^{\prime}=(0,0,0,0,\frak d).$$
That $\varphi$ is an isomorphism is verified by showing that 
$\{e_1^{\prime},e_2^{\prime},e_3^{\prime},e_4^{\prime},e_5^{\prime}
\}\subset \roman H_{5,3}$ satisfy $(C)$. (Here $\varphi$ is given by 
$$(h,j,k,m,n)\mapsto (\delta \frak d^2h+(\gamma\delta \frak d(\frak d
-1)/2)j,\gamma\delta \frak dj+\delta\epsilon \frak dk+\beta\delta 
\frak dm,\delta\epsilon \frak dk,\alpha\gamma\delta m,\frak dn).)$$
When $\epsilon =0$, use $\frak d=\alpha\gamma$ and $e_3^{\prime}=(
0,0,\delta \frak d,0,0)$. 

\smallpagebreak

It is easy to see that the image $\roman H_1=\varphi (\roman H_{5,
3}(\alpha ,\beta ,\gamma ,\delta ,\epsilon ))$ is cofinite in 
$\roman H_{5,3}$. Consider the coset $s\roman H_1$ for $s=(h,j,k,m
,n)\in \roman H_{5,3}$; since 
$e_5^{\prime}=(0,0,0,0,\frak d)$, we can choose 
$r_5\in \Bbb Z$ so that $se^{\prime}_5{}^{r_5}$ has its last coordinate in $
[0,\frak d)$. Then choose $r_4\in \Bbb Z$ 
so that $se^{\prime}_5{}^{r_5}e^{\prime}_4{}^{r_4}$ has its second last coordinate in $
[0,\alpha\gamma\delta )$. Continuing 
like this, we arrive at 
$$se^{\prime}_5{}^{r_5}e^{\prime}_4{}^{r_4}e^{\prime}_3{}^{r_3}e^{
\prime}_2{}^{r_2}e^{\prime}_1{}^{r_1}\in$$
$$K=\left([0,\delta \frak d^2)\times [0,\gamma\delta \frak d)\times 
[0,\delta\epsilon \frak d)\times [0,\alpha\gamma\delta )\times [0,
\frak d)\right)\cap \Bbb Z^5\subset \roman H_{5,3},$$
so every coset $s\roman H_1$, $s\in \roman H_{5,3}$, has a representative 
in $K$, which is a finite set. It follows that the quotient map 
$\roman H_{5,3}\to \roman H_{5,3}/\roman H_1$ maps $K$ onto 
$\roman H_{5,3}/\roman H_1$, which is therefore finite. (A similar argument shows that 
$\roman G_{5,3}/\roman H_1$ is cocompact.) 

Finally, note that since any cofinite subgroup of $\roman H_{5,3}$ is also a discrete 
cocompact subgroup of $\roman G_{5,3}$, it must therefore be isomorphic to some 
$\roman H_{5,3}(\alpha ,\beta ,\gamma ,\delta ,\epsilon )$. This completes the proof. \qed

\bigpagebreak

\flushpar REMARKS. 1. The image $\roman H_1=\varphi (\roman H_{5,3}
(\alpha ,\beta ,\gamma ,\delta ,\epsilon ))$ above is not a normal 
subgroup of $\roman H_{5,3}$, e.g.,
$$(0,0,1,0,0)e_5^{\prime}(0,0,-1,0,0)=(\frak d,0,0,0,0)\notin \roman H_
1.$$
This makes it seem unlikely that $\roman H_{5,3}(\alpha ,\beta ,\gamma 
,\delta ,\epsilon )$ can be embedded in 
$\roman H_{5,3}$ as a normal subgroup; however, the existence of such an 
embedding is still a possiblity.

\smallpagebreak

\flushpar 2. The theorem gives an isomorphism $\varphi$ of
$\roman H_{5,3}(\alpha ,\beta ,\gamma ,\delta ,\epsilon )$ into $\roman H_{
5,3}$; conversely, there is always
an isomorphism $\varphi^{\prime}$ of $\roman H_{5,3}$ into $\roman H_{
5,3}(\alpha ,\beta ,\gamma ,\delta ,\epsilon )$, and as for 
$\varphi$, it is easier to give $\varphi^{\prime}$ in terms of the generators $
e_i$, 
$1\leq i\leq 5$, of $\roman H_{5,3}$, which satisfy 
$$[e_5,e_4]=e_2,\,\,\,[e_5,e_2]=e_1=[e_4,e_3].\leqno (C^{\prime})$$
Then
$$\varphi^{\prime}:e_1\mapsto e_1^{\prime}=(\alpha\gamma^2\delta^2
,0,0,0,0),\,\,\,e_2\mapsto e_2^{\prime}=(\alpha\gamma^2\delta (\delta 
-1)/2,\alpha\gamma\delta ,0,0,0),$$
$$e_3\mapsto e_3^{\prime}=(0,-\alpha\delta\epsilon ,\alpha\delta\gamma 
,0,0),\,\,\,e_4\mapsto e_4^{\prime}=(0,-\beta ,0,\gamma ,0),$$
$$\text{ and }e_5\mapsto e_5^{\prime}=(0,0,0,0,\delta ).$$
That $\varphi^{\prime}$ is an isomorphism  is verified by showing that 
$\{e_1^{\prime},e_2^{\prime},e_3^{\prime},e_4^{\prime},e_5^{\prime}
\}\subset \roman H_{5,3}(\alpha ,\beta ,\gamma ,\delta ,\epsilon )$ satisfy $
(C^{\prime})$.
(Here $\varphi^{\prime}$ is given by 
$$(h,j,k,m,n)\mapsto (\alpha\gamma^2\delta^2h+j\alpha\gamma^2\delta 
(\delta -1)/2,\alpha\gamma\delta j-\alpha\delta\epsilon k-\beta m,
\alpha\gamma\delta k,\gamma m,\delta n).)$$

So, as for the 3-dimensional groups $\roman H_3(p)$
and the 4-dimensional groups $\roman H_4(p_1,p_2,p_3)$, here we have an infinite 
family of non-isomorphic groups, each of
which is isomorphic to a subgroup of any other one.


\subhead \S 3. Infinite Dimensional  Simple Quotients of 
$C^{*}(\roman H_{5,3}(\alpha ,\beta ,\gamma ,\delta ,\epsilon ))$
\endsubhead

We begin by obtaining concrete representations on $L^2(\Bbb T^2)$ of the faithful simple 
quotients (i.e., those arising from a faithful representation of 
$\roman H_{5,3}(\alpha ,\beta ,\gamma ,\delta ,\epsilon )$), and consider first the case $
\epsilon =0$. In this case 
$\roman H_{5,3}(\alpha ,\beta ,\gamma ,\delta ,0)$ has an 
abelian normal subgroup $N=(\Bbb Z,\Bbb Z,0,\Bbb Z,0)$, with quotient 
$$\roman H_{5,3}(\alpha ,\beta ,\gamma ,\delta ,0)/N\cong (0,0,\Bbb Z
,0,\Bbb Z)=\Bbb Z^2,$$
also abelian and embedded in $\roman H_{5,3}(\alpha ,\beta ,\gamma 
,\delta ,0)$ as a subgroup, so that 
$\roman H_{5,3}(\alpha ,\beta ,\gamma ,\delta ,0)$ is isomorphic to a semidirect product $
N\times \Bbb Z^2$; in this 
situation, the simple quotients of $C^{*}(\roman H_{5,3}(\alpha ,\beta 
,\gamma ,\delta ,0))$ can be 
presented as $C^{*}$-crossed products using flows from commuting 
homeomorphisms, as follows.

\smallpagebreak

Let $\lambda =e^{2\pi i\theta}$ for an irrational $\theta$, and consider 
the flow $\Cal F^{\prime}=(\Bbb Z^2,\Bbb T^2)$ generated by the commuting homeomorphisms
$$\psi^{\prime}_1:(w,v)\mapsto (\lambda^{\gamma}w,\lambda^{\beta}w^{
\alpha}v)\,\,\,\text{ and }\,\,\,\psi_2^{\prime}:(w,v)\mapsto (w,\lambda^{
-\delta}v).$$
$\Cal F^{\prime}$ is minimal, so the $C^{*}$-crossed product $\Cal C^{
\prime}=C^{*}(\Cal C(\Bbb T^2),\Bbb Z^2)$ is 
simple [\EH , Corollary 5.16]. 

Let $v$ and $w$ denote (as well as members of $\Bbb T$) the functions in 
$\Cal C(\Bbb T^2)$ defined by
$$(w,v)\mapsto v\text{ and }w,$$
respectively. Define unitaries $U$, $V$, $W$ and $X$ on $L^2(\Bbb T^
2)$ by 
$$U:f\mapsto f\circ\psi_1^{\prime},\,\,\,\,\,\,\,V:f\mapsto vf,\,\,\,\,\,\,
W:f\mapsto f\circ\psi^{\prime}_2\,\,\,\text{ and }\,\,\,X:f\mapsto 
wf.\leqno (\Cal U^{\prime})$$
These unitaries satisfy 
$$UV=\lambda^{\beta}X^{\alpha}VU,\,\,\,\,\,\,\,UX=\lambda^{\gamma}
XU,\,\,\,\text{ and }\,\,\,\,VW=\lambda^{\delta}WV\leqno (\roman C
\roman R^{\prime})$$
(other pairs of unitaries commuting), equations which ensure that 
$$\pi :(h,j,k,m,n)\mapsto\lambda^hX^jW^kV^mU^n$$
is a representation of $\roman H_{5,3}(\alpha ,\beta ,\gamma ,\delta 
,0)$. Denote by $\roman A_{\theta}^{5,3}(\alpha ,\beta ,\gamma ,\delta 
,0)$ the  
$C^{*}$-subalgebra of $B(L^2(\Bbb T^2))$ generated by $\pi$, i.e., by $
U$, $V$, $W$ and $X$.
Since $\roman A_{\theta}^{5,3}(\alpha ,\beta ,\gamma ,\delta ,0)$ is generated by a representation of 
$\roman H_{5,3}(\alpha ,\beta ,\gamma ,\delta ,0)$, it 
is a quotient of the group $C^{*}$-algebra 
$C^{*}(\roman H_{5,3}(\alpha ,\beta ,\gamma ,\delta ,0))$. It follows readily that $
\roman A_{\theta}^{5,3}(\alpha ,\beta ,\gamma ,\delta ,0)$ is isomorphic 
to the simple $C^{*}$-crossed product $\Cal C^{\prime}$ above, and hence is simple.

\smallpagebreak

However, when $0<\epsilon\leq \roman g\roman c\roman d\,\{\gamma ,
\delta \}/2$ (which implies $\gamma >1$, by $(*)$), 
$\roman H_{5,3}(\alpha ,\beta ,\gamma ,\delta ,\epsilon )$ is only an extension
$(\Bbb Z,\Bbb Z,0,\Bbb Z,0)\times (0,0,\Bbb Z,0,\Bbb Z)=N\times \Bbb Z^
2$, and not a semidirect product. 
Nonetheless, we can modify the flow $\Cal F^{\prime}$ representing $
\roman A_{\theta}^{5,3}(\alpha ,\beta ,\gamma ,\delta ,0)$
above to get a concrete representation of $\roman A_{\theta}^{5,3}
(\alpha ,\beta ,\gamma ,\delta ,\epsilon )$.
Consider the flow $\Cal F=(\Bbb Z^2,\Bbb T^2)$ generated by the commuting homeomorphisms
$$\psi_1:(w,v)\mapsto (\lambda w,\lambda^{\beta}w^{\alpha\gamma}v)\text{ and }
\psi_2:(w,v)\mapsto (w,\lambda^{-\delta}v).$$
$\Cal F$ is minimal, so the $C^{*}$-crossed product $\Cal C=C^{*}(
\Cal C(\Bbb T^2),\Bbb Z^2)$ is simple.
Define unitaries on $L^2(\Bbb T^2)$ by
$$U:f\mapsto f\circ\psi_1,\,\,\,\,\,\,V:f\mapsto vf,\,\,\,\,\,\,W:
f\mapsto w^{\epsilon}f\circ\psi_2\,\,\,\text{ and }\,\,\,\,X:f\mapsto 
w^{\gamma}f.\leqno (\Cal U)$$
These unitaries satisfy 
$$UV=\lambda^{\beta}X^{\alpha}VU,\,\,\,\,\,\,UX=\lambda^{\gamma}XU
,\,\,\,\,\,\,VW=\lambda^{\delta}WV\,\,\,\text{ and }\,\,\,\,UW=\lambda^{
\epsilon}WU,\leqno (\roman C\roman R)$$
equations which ensure that 
$$\pi :(h,j,k,m,n)\mapsto\lambda^hX^jW^kV^mU^n$$
is a representation of $\roman H_{5,3}(\alpha ,\beta ,\gamma ,\delta 
,\epsilon )$. Denote by $\roman A_{\theta}^{5,3}(\alpha ,\beta ,\gamma 
,\delta ,\epsilon )$ the  
$C^{*}$-subalgebra of $B(L^2(\Bbb T^2))$ generated by $\pi$. Now $
\roman A_{\theta}^{5,3}(\alpha ,\beta ,\gamma ,\delta ,\epsilon )$ is
isomorphic only to a subalgebra of $\Cal C$ (as may be shown using 
conditional expectations); a unitary that is 
missing is $X^{\prime}:f\mapsto wf$ (since $\gamma >1)$. 

\smallpagebreak

\flushpar NOTE. The reason we did not use $\Cal F$ when $\epsilon 
=0$ (and $\gamma >1)$ is that 
$\roman A_{\theta}^{5,3}(\alpha ,\beta ,\gamma ,\delta ,0)$ seems to be isomorphic 
only to a subalgebra of $\Cal C$ in that case too, whereas with 
$\Cal F^{\prime}$, 
$\roman A_{\theta}^{5,3}(\alpha ,\beta ,\gamma ,\delta ,0)\cong \Cal C^{
\prime}$.

\smallpagebreak

Since the flow method can no longer be used to prove 
the simplicity of the algebra 
$\roman A_{\theta}^{5,3}(\alpha ,\beta ,\gamma ,\delta ,\epsilon )$ 
(when $0<\epsilon\leq \roman g\roman c\roman d\,\{\gamma ,\delta \}
/2$),
we use the strong result of Packer [\JPa ].

\proclaim{2.~Theorem} Let $\lambda =e^{2\pi i\theta}$ for an irrational $
\theta$.

\flushpar$(a)$ There is a unique $($up to isomorphism$)$ simple
$C^{*}$-algebra $\roman A_{\theta}^{5,3}(\alpha ,\beta ,\gamma ,\delta 
,\epsilon )$ generated by unitaries $U$, $V$, $W$ 
and $X$ satisfying 
$$UV=\lambda^{\beta}X^{\alpha}VU,\,\,\,\,\,\,UX=\lambda^{\gamma}XU
,\,\,\,\,\,\,VW=\lambda^{\delta}WV\,\,\,\text{ and }\,\,\,\,UW=\lambda^{
\epsilon}WU,\leqno (\roman C\roman R)$$
Furthermore, for a suitable $\Bbb C$-valued cocycle on $\roman H_3
(\alpha )\times \Bbb Z$, 
$$\roman A_{\theta}^{5,3}(\alpha ,\beta ,\gamma ,\delta ,\epsilon 
)\cong C^{*}(\Bbb C,\roman H_3(\alpha )\times \Bbb Z).$$

\flushpar$(b)$ Let $\pi^{\prime}$ be a representation of 
$\roman H_{5,3}^{\prime}=\roman H_{5,3}(\alpha ,\beta ,\gamma ,\delta 
,\epsilon )$ such that $\pi =\pi^{\prime}$ $($as
scalars$)$ on the center $(\Bbb Z,0,0,0,0)$ of $\roman H_{5,3}^{\prime}$, and let 
$\roman A$ be the $C^{*}$-algebra
generated by $\pi^{\prime}$. Then $\roman A\cong \roman A_{\theta}^{
5,3}(\alpha ,\beta ,\gamma ,\delta ,\epsilon )=\roman A^{\prime 5,
3}_{\theta}$ $($say$)$ via a unique 
isomorphism $\omega$ such that the following diagram commutes.
$$\matrix \roman H_{5,3}^{\prime}\,\,\longrightarrow\!\!\!\!\!\!\!\!^{
\pi}&\,\,\roman A^{\prime 5,3}_{\theta}\\
&\!\!\!\!\!\!\!\!\!\!\!\!\!\!\!\!\!\!\!\!\!\!\!\!\!\!\!\!\!\!\!\!\!\!\!\!\!\!\!\!\!\!\!\!\!\!\!\!\!\!\!\!
\pi^{\prime}\searrow&\!\!\!\!\!\!\!\!\!\!\!\!\!\!\!\!\!\!\!\!\!\!\!\!
\swarrow\omega\\
&\!\!\!\!\!\!\!\!\!\!\!\!\!\!\!\!\!\!\!\!\!\!\!\!\!\!\!\!\!\roman A\endmatrix$$
\endproclaim

\flushpar PROOF. 
To use Packer's result, we regard $\roman H_{5,3}(\alpha ,\beta ,\gamma 
,\delta ,\epsilon )$ as an 
extension 
$$\roman H_{5,3}(\alpha ,\beta ,\gamma ,\delta ,\epsilon )\cong \Bbb Z
\times (0,\Bbb Z,\Bbb Z,\Bbb Z,\Bbb Z)\cong \Bbb Z\times (\roman H_
3(\alpha )\times \Bbb Z)$$
(with $\roman H_3(\alpha )\cong (0,\Bbb Z,0,\Bbb Z,\Bbb Z)\subset 
\roman H_{5,3}(\alpha ,\beta ,\gamma ,\delta ,\epsilon )$); this extension has cocycle 
$$[s,s^{\prime}]=[(j,k,m,n),(j^{\prime},k^{\prime},m^{\prime},n^{\prime}
)]=\lambda^{\gamma nj^{\prime}+\alpha\gamma m^{\prime}n(n-1)/2+\beta 
nm^{\prime}+\delta mk^{\prime}+\epsilon nk^{\prime}},$$
$$(\roman H_3(\alpha )\times \Bbb Z,\roman H_3(\alpha )\times \Bbb Z
)\to \Bbb T.$$
The application of Packer's result requires the consideration of the related function  
$$\chi^{s^{\prime}}(s)=[s^{\prime},s]\overline {[s,s^{-1}s^{\prime}
s]},\,\,\,s,s^{\prime}\in (0,\Bbb Z,\Bbb Z,\Bbb Z,\Bbb Z)\cong \roman H_
3(\alpha )\times \Bbb Z.$$
It must be shown that $\chi^{s^{\prime}}$ is non-trivial on the centralizer of $
s^{\prime}$ in $\roman H_3(\alpha )\times \Bbb Z$ if 
$s^{\prime}$ has finite conjugacy class in $\roman H_3(\alpha )\times 
\Bbb Z$; this is easy because the 
only elements of $\roman H_3(\alpha )\times \Bbb Z$ that have finite conjugacy class are in the 
center $Z_1=(\Bbb Z,\Bbb Z,0,0)$ of $\roman H_3(\alpha )\times \Bbb Z$, so their centralizer is all of 
$\roman H_3(\alpha )\times \Bbb Z$. Thus the $C^{*}$-crossed product $
C^{*}(\Bbb C,\roman H_3(\alpha )\times \Bbb Z)$ is simple; it is 
isomorphic to $\roman A_{\theta}^{5,3}(\alpha ,\beta ,\gamma ,\delta 
,\epsilon )$ because, with basis members
$$e_1=(1,0,0,0),\,\,\,\,\,\,e_2=(0,1,0,0),\,\,\,\,\,\,e_3=(0,0,1,0
)\,\,\,\text{ and }\,\,\,\,e_4=(0,0,0,1)$$
for $\roman H_3(\alpha )\times \Bbb Z$, the unitaries 
$$U^{\prime}=\delta_{e_4},\,\,\,\,\,\,V^{\prime}=\delta_{e_3},\,\,\,\,\,\,
W^{\prime}=\delta_{e_2}\,\,\,\text{ and }\,\,\,\,X^{\prime}=\delta_{
e_1}$$
in $\ell_1(\roman H_3(\alpha )\times \Bbb Z)\subset C^{*}(\Bbb C,\roman H_
3(\alpha )\times \Bbb Z)$ satisfy (CR). \qed

\bigpagebreak

The theorem showed $\roman A_{\theta}^{5,3}(\alpha ,\beta ,\gamma 
,\delta ,\epsilon )$ was simple by showing it was 
isomorphic to the simple $C^{*}$-crossed product $C^{*}(\Bbb C,\roman H_
3(\alpha )\times \Bbb Z)$. It 
follows that $\roman A_{\theta}^{5,3}(\alpha ,\beta ,\gamma ,\delta 
,\epsilon )$ is isomorphic to a number of other 
$C^{*}$-crossed products (much as in [\AM ; Theorem 3]), one of which has been 
derived from a flow at the beginning of the section for the case 
$\epsilon =0$. Here are 2 other $C^{*}$-crossed products that can be used to
establish the simplicity of $\roman A_{\theta}^{5,3}(\alpha ,\beta 
,\gamma ,\delta ,\epsilon )$. The variable change for the 
second one will be used again, in the proof of Theorem 3 below.

\smallpagebreak

\flushpar 1. Take the $C^{*}$-algebra $B$ generated by $U$, $V$ and $
X$ from (CR), satisfying 
$$UV=\lambda^{\beta}X^{\alpha}VU,\,\,\,\,\,\,UX=\lambda^{\gamma}XU
,\,\,\,\text{ and }\,\,\,\,VX=XV.$$
The algbera $B$ is a faithful simple quotient of a discrete cocompact subgroup 
$\roman H_4(\beta ,\alpha ,\gamma )$ of $\roman H_4$, the connected 4-dimensional nilpotent group 
[\MWc ; Theorem 2].  Then the rest of (CR) gives an action of $\Bbb Z$ on 
$B$ generated by $\roman A\roman d_W$; it follows that $\roman A_{
\theta}^{5,3}(\alpha ,\beta ,\gamma ,\delta ,\epsilon )\cong C^{*}
(B,\Bbb Z)$. The 
simplicity of $C^{*}(B,\Bbb Z)$ can be proved directly by showing that 
$(\roman A\roman d_W)^r=\text{$\roman A\roman d_{W^r}$}$ is outer on $
B$ if $r\neq 0$ [8].

\smallpagebreak

\flushpar 2. 
First, we change the variables in (CR).  Pick relatively prime integers $
c,d$ such that 
$d\delta +c\epsilon =0$ and let $a,b$ be integers such that $ad-bc
=1$.  Put
$$U^{\prime}=U^aV^b\,\,\,\text{ and }\,\,\,V^{\prime}=U^cV^d.$$
Then keeping $X$ and $W$ the same, $(\roman C\roman R$) becomes
$$\left\{\aligned
&U^{\prime}V^{\prime}=\lambda^{\beta^{\prime}}X^{\alpha}V^{\prime}
U^{\prime},\,\,\,\,\,\,U^{\prime}X=\lambda^{a\gamma}XU^{\prime},\\
&U^{\prime}W=\lambda^{\delta^{\prime}}WU^{\prime},\,\,\,\text{ and }\,\,\,
V^{\prime}X=\lambda^{c\gamma}XV^{\prime}\endaligned
\right.\leqno (\roman C\roman R^{\prime})$$
(other pairs of unitaries commuting) for some integer $\beta^{\prime}$ and 
$\delta^{\prime}=b\delta +a\epsilon$.  Note that $\roman A_{\theta}^{
5,3}(\alpha ,\beta ,\gamma ,\delta ,\epsilon )$ is generated by $U^{
\prime},V^{\prime},W,X$ (since 
$ad-bc=1$). Consider the $C^{*}$-algebra $B^{\prime}=\roman A_{\theta^{
\prime}}\otimes \Cal C(\Bbb T)$ generated by 
unitaries $V^{\prime}$, $W$ and $X_1$ satisfying
$$V^{\prime}X=\lambda^{c\gamma}XV^{\prime},\,\,\,\,\,\,V^{\prime}W
=WV^{\prime}\,\,\,\text{ and }\,\,\,\,WX_1=X_1W;$$
here $e^{2\pi i\theta^{\prime}}=\lambda^{c\gamma}$.
Then the rest of $(\roman C\roman R^{\prime})$ gives an action of $
\Bbb Z$ on $B^{\prime}$ generated by $\roman A\roman d_{U^{\prime}}$, 
and it follows 
that $C^{*}(B^{\prime},\Bbb Z)\cong \roman A_{\theta}^{5,3}(\alpha 
,\beta ,\gamma ,\delta ,\epsilon )$. One can prove the simplicity of 
$C^{*}(B^{\prime},\Bbb Z)$ directly; the method of proof is to show that the Connes 
spectrum of $\roman A\roman d_{U^{\prime}}$ is $\Bbb T$, which follows from Theorem 2 and 
[\GP ; 8.11.12].

\bigpagebreak


\subhead \S 4. Other Simple Quotients of $C^{*}(\roman H_{5,3}(\alpha 
,\beta ,\gamma ,\delta ,\epsilon ))$
\endsubhead

Now assume that $\lambda$ is a primitive $q$th root of unity and that
$U$, $V$, $W$ and $X$ are unitaries generating a simple quotient $
A$ 
of $C^{*}(\roman H_{5,3}(\alpha ,\beta ,\gamma ,\delta ,\epsilon )
)$, i.e., they satisfy 
$$UV=\lambda^{\beta}X^{\alpha}VU,\,\,\,\,\,\,UX=\lambda^{\gamma}XU
,\,\,\,\,\,\,VW=\lambda^{\delta}WV\,\,\,\text{ and }\,\,\,\,UW=\lambda^{
\epsilon}WU.\leqno (\roman C\roman R)$$
We may assume that $A$ is irreducibly represented. 
Then, if 
$$\left\{\aligned
&q_1\text{ is the order of }\lambda^{\gamma}\text{ and}\\
&q_2\text{ is the lcm of the orders of }\lambda^{\delta}\text{ and }
\lambda^{\epsilon},\endaligned
\right.\leqno (c^{\prime})$$
$W^{q_2}$ and $X^{q_1}$ are scalar multiples of the identity (by irreducibility). 
Since $W$ can be multiplied by a scalar without changing (CR), we may 
assume $W^{q_2}=1$. However, $X^{q_1}=\mu^{\prime}$, a multiple of the identity. 
Put $X=\mu X_1$ for $\mu^{q_1}=\mu^{\prime}$, so that $X_1^{q_1}=1$, and
substitute $X=\mu X_1$ in (CR) to get
$$\left\{\aligned
&UV=\lambda^{\beta}\mu^{\alpha}X_1^{\alpha}VU,\,\,\,\,\,\,UX_1=\lambda^{
\gamma}X_1U,\\
&VW=\lambda^{\delta}WV,\,\,\,\,\,\,UW=\lambda^{\epsilon}WU\,\,\,\text{ and }\,\,\,\,
W^{q_2}=1=X^{q_1}_1.\endaligned
\right.\leqno (\roman C\roman R_1)$$

\smallpagebreak

1. If $\mu$ is also a root of unity, then $(\roman C\roman R_1)$ (along with irreducibility) 
shows that $U$ and $V$, as well 
as $W$ and $X$, are (multiples of) unipotent unitaries, so $A$ is finite dimensional.

\smallpagebreak

2. If $\mu$ is not a root of unity,
the flow $\Cal F=(\Bbb Z^2,\Bbb T^2)$ used above to get a concrete representation of 
$\roman A_{\theta}^{5,3}(\alpha ,\beta ,\gamma ,\delta ,\epsilon )$ can be modified to get a concrete representation of $
A$ 
on $L^2(\Bbb Z_{q_1}\times \Bbb T)$ (where $\Bbb Z_{q_1}$ is the subgroup of $
\Bbb T$ with $q_1$ elements). The 
proof of the simplicity of $A$ comes next.

\bigpagebreak

First consider the universal $C^{*}$-algebra $\frak A$ generated by unitaries satisfying
$$\left\{\aligned
&UV=\lambda^{\beta}\mu^{\alpha}X_1^{\alpha}VU,\,\,\,\,\,\,UX_1=\lambda^{
\gamma}X_1U,\\
&VW=\lambda^{\delta}WV,\,\,\,\,\,\,UW=\lambda^{\epsilon}WU\,\,\,\text{ and }\,\,\,\,
W^{q_2}=1=X^{q_1}_1.\endaligned
\right.\leqno (\roman C\roman R_1)$$
A change of variables is useful.
Pick relatively prime integers $c,d$ such that
$d\delta +c\epsilon =0$ and let $a,b$ be integers such that $ad-bc
=1$.  Put
$$U^{\prime}=U^aV^b\,\,\,\text{ and }\,\,\,\,V^{\prime}=U^cV^d.$$
Then keeping $X$ and $W$ the same, (CR${}_1$) 
becomes
$$\left\{\aligned
&U^{\prime}V^{\prime}=\xi X_1^{\alpha}V^{\prime}U^{\prime},\,\,\,\,\,\,
U^{\prime}X_1=\lambda^{a\gamma}X_1U^{\prime},\\
&U^{\prime}W=\lambda^{\delta^{\prime}}WU^{\prime},\,\,\,\,\,\,V^{\prime}
X=\lambda^{c\gamma}XV^{\prime}\,\,\,\text{ and }\,\,\,\,W^{q_2}=1=
X^{q_1}_1\endaligned
\right.\leqno (\roman C\roman R_2)$$
(other pairs of unitaries commuting),
where $\xi =\lambda^{\beta}\mu^{\alpha}\lambda^s$ for some integer $
s$, and 
$\delta^{\prime}=b\delta +a\epsilon$. It is clear that $\lambda^{\delta^{
\prime}}$ is a primitive
$q_2$-th root of unity and
that the algebra $\frak A$ is generated by $U^{\prime},V^{\prime},
W$ and $X_1$, since $ad-bc=1$.

Let $B=C^{*}(X_1,V^{\prime})$ and let $C(\Bbb Z_{q_2})=C^{*}(W)$ be the $
C^{*}$-algebra generated by $W$.
Since $W$ commutes with $X_1$ and $V^{\prime}$, we can form the tensor product algebra
$B\otimes C(\Bbb Z_{q_2})=C^{*}(X_1,V^{\prime},W)$.  The automorphism $
\roman A\roman d_{U^{\prime}}$ acts on this
tensor product as $\sigma\otimes\tau$, where $\sigma$ and $\tau$ are automorphisms of
$B$ and $C(\Bbb Z_{q_2})$, respectively, given by
$$\sigma (X_1)=\lambda^{a\gamma}X_1,\,\,\,\,\,\,\sigma (V^{\prime}
)=\xi_1X_1^{\alpha}V^{\prime}\,\,\,\text{ and }\,\,\,\,\tau (W)=\zeta 
W.$$
Therefore, by the universality of $\frak A$ and of the $C^{*}$-crossed product 
$C^{*}(B\otimes C(\Bbb Z_{q_2}),\Bbb Z)$, these algebras are isomorphic.
By Rieffel's Proposition 1.2 [\MR ], 
the latter of these is isomorphic to $M_{q_2}(D)$, where
$D=C^{*}(B,\Bbb Z)=C^{*}(X_1,V^{\prime},U^{\prime}{}^{q_2})$, and the action of $
\Bbb Z$ on $B$ is generated by 
$\sigma^{q_2}$.

Now, the unitaries $X_1,V^{\prime}\text{ and }U^{\prime}{}^{q_2}$ generating $
D$ satisfy
$$\left\{\aligned
&U^{\prime}{}^{q_2}V^{\prime}=\xi^{q_2}\lambda^{s^{\prime}}X_1^{\alpha 
q_2}V^{\prime}U^{\prime}{}^{q_2},\,\,\,\,\,\,V^{\prime}X_1=\lambda^{
c\gamma}X_1V^{\prime},\\
&U^{\prime}{}^{q_2}X_1=\lambda^{a\gamma q_2}X_1U^{\prime}{}^{q_2}\,\,\,\text{ and }\,\,\,\,
X_1^{q_1}=1,\endaligned
\right.\leqno (\star )$$
for some $s^{\prime}\in \Bbb Z$.

Now we apply another change of
variables.  Choose relatively prime integers $c^{\prime},d^{\prime}$ such that 
$cd^{\prime}+aq_2c^{\prime}=0$, then pick integers $a^{\prime},b^{
\prime}$ with $a^{\prime}d^{\prime}-b^{\prime}c^{\prime}=1$, and put
$$U^{\prime\prime}=U^{\prime}{}^{q_2a^{\prime}}V^{\prime}{}^{b^{\prime}}\,\,\,\text{ and }\,\,\,\,
V^{\prime\prime}=U^{\prime}{}^{q_2c^{\prime}}V^{d^{\prime}}.$$
Then $(\star$) becomes (keeping $X_1$ the same)
$$\left\{\aligned
&U^{\prime\prime}V^{\prime\prime}=\xi_1X_1^{\alpha q_2}V^{\prime\prime}
U^{\prime\prime},\,\,\,\,\,\,\,\,\,V^{\prime\prime}X_1=X_1V^{\prime
\prime},\\
&U^{\prime\prime}X_1=\lambda^{\prime}X_1U^{\prime\prime}\,\,\,\text{ and }\,\,\,\,
X_1^{q_1}=1,\endaligned
\right.\leqno (\star\star )$$
where $\xi_1=\xi^{q_2}\lambda^{s^{\prime}}$ for some integer $s^{\prime}$,
$\lambda^{\prime}=\lambda^{\gamma (aq_2a^{\prime}+cb^{\prime})}$ has order $
q_3$ dividing $q_1$ (the order of $\lambda^{\gamma}$), and 
perhaps $q_3\neq q_1.$

Now, with $\Bbb Z_{q_1}\subset \Bbb T$ representing the subgroup with $
q_1$ members,  one 
observes that $D$ is isomorphic to the crossed product of 
$C^{*}(C(\Bbb Z_{q_1}\times \Bbb T),\Bbb Z)$ from the flow generated by
$\phi (w,v)=(\lambda^{\prime}w,\xi_1\lambda^{\gamma\alpha q_2}v)$.  (Note that the flow is not minimal
unless the order of $\lambda^{\prime}$ is exactly $q_1$.) 
This proves the following. 

\proclaim{3. Theorem} The universal $C^{*}$-algebra $\frak A$ generated by
unitaries $U$, $V$, $W$ and $X_1$ satisfying $(\roman C\roman R_1)$
as for $2$ $($see also $($c$^{\prime}))$ is isomorphic to $M_{q_2}
(D)$, where 
$D=C^{*}(C(\Bbb Z_{q_1}\times \Bbb T),\Bbb Z)$, as above. \endproclaim

Therefore, we now obtain all simple algebras satisfying $(\roman C
\roman R_1)$.

\proclaim{4. Corollary}
Every simple $C^{*}$-algebra generated by unitaries satisfying $(\roman C
\roman R_1)$ is isomorphic to a
matrix algebra over an irrational rotation algebra.
\endproclaim

\flushpar PROOF.
By Theorem 3, any such simple algebra $Q$ is a quotient of $M_{q_2}
(D)$.  
Hence $Q=M_{q_2}(Q^{\prime})$ where $Q^{\prime}$ is a simple quotient of $
D$.  But such a
$Q^{\prime}$ is generated by unitaries satisfying $(\star\star$), but with $
X_1$ (of order $q_1$) replaced by 
another unitary $X_2$, which
after suitable rescaling, has order equal to the order of $\lambda^{
\prime}$. But this 
algebra is known to be a matrix algebra over an 
irrational rotation algebra
(see for example Theorem 3 of [\MWa ]). \qed

\bigpagebreak

We state

\proclaim{5. Theorem} A $C^{*}$-algebra $\roman A$ is isomorphic to a 
simple infinite dimensional quotient of $C^{*}(\roman H_{5,3}(\alpha 
,\beta ,\gamma ,\delta ,\epsilon ))$
if and only if $\roman A$ is isomorphic to $\roman A_{\theta}^{5,3}
(\alpha ,\beta ,\gamma ,\delta ,\epsilon )$ for an irrational 
$\theta$, or to an algebra as in Corollary $4$.
\endproclaim

\bigpagebreak

The matrix algebra presentation for the simple $C^{*}$-algebra $A$ generated by 
unitaries satisfying $(\roman C\roman R_1)$ is not given in a definite form in 
the proof of Corollary 4. Here is an explicit matrix presentation. 
First, as near the beginning of the section, change the 
variables so that $(\roman C\roman R_1)$ becomes
$$\left\{\aligned
&U^{\prime}V^{\prime}=\xi X_1^{\alpha}V^{\prime}U^{\prime},\,\,\,\,\,\,
U^{\prime}X_1=\lambda^{a\gamma}X_1U^{\prime},\\
&U^{\prime}W=\lambda^{\delta^{\prime}}WU^{\prime},\,\,\,\,\,\,V^{\prime}
X=\lambda^{c\gamma}XV^{\prime}\,\,\,\text{ and }\,\,\,\,W^{q_2}=1=
X^{q_1}_1\endaligned
\right.\leqno (\roman C\roman R_2)$$

Now we shall present the algebra $A$ by unitaries in a matrix algebra as follows.
Consider the $C^{*}$-algebra $B_1$ generated by unitaries 
$u,v$ and $x$ enjoying the relations
$$uv=\xi^{\prime}x^{q_2\alpha}vu,\,\,\,ux=\lambda^{q_2a\gamma}xu,\,\,\,
vx=\lambda^{c\gamma}xv\text{ and }x^{q_3}=1,$$
where $q_3$ is the least common multiple of the orders of $\lambda^{
q_2a\gamma}$ and $\lambda^{c\gamma}$, 
and $\xi^{\prime}$ is to be determined below.  Clearly, $q_3$ divides $
q_1$ so that also
$x^{q_1}=1$.  It was shown in the proof of Theorem 6.4 of [\MWc ] 
that $B_1$ is 
isomorphic to a $q_3\times q_3$ matrix algebra over an irrational rotation algebra,
when $\xi^{\prime}$ is not a root of unity.  Hence it will suffice to show that $
A$ is 
isomorphic to $M_{q_2}(B)$ (so that $Q=q_2q_3$).  Indeed, let
$$
\split
U^{\prime} &= 
\left(\matrix 0&0&\cdots&0&u\\
1&0&\cdots&0&0\\
0&1&\cdots&0&0\\
\vdots&&\ddots&&\vdots\\
0&0&\cdots&1&0\endmatrix\right)
\intertext{(so $U^{\prime}\text{diag}(K_1,K_2,\dots ,K_{q_2})U^{\prime}{}^{*}
=\text{diag}(uK_{q_2}u^{*},K_1,K_2,\dots ,K_{q_2-1}))$,}
V^{\prime} &=\text{diag}(\tau_1v,~\tau_2x^{-\alpha}v,~\tau_3x^{-2\alpha}
v,\dots ,~\tau_{q_2}x^{-(q_2-1)\alpha}v),
\\
W &= \text{diag}(1,\zeta^{-1},\zeta^{-2},\dots ,\zeta^{-(q_2-1)}),
\\
X &= \text{diag}(x,~\lambda^{-a\gamma}x,~\lambda^{-2a\gamma}x,\dots,
~\lambda^{-(q_2-1)a\gamma}x),
\endsplit
$$ 
where 
$$\tau_j=\lambda^{a\alpha\gamma j(j-1)/2}\,\xi_1^{-j+1},\,\,\,1\leq 
j\leq q_2,\text{ and }\tau_{q_2}\xi^{\prime}\lambda^{-q_2(q_2-1)\alpha^
2\gamma}=\xi_1.$$
One now checks that these unitaries satisfy (CR${}_2$). It is also 
evident that they generate $M_{q_2}(B)$.  \qed

\bigpagebreak

\newpage

\subhead \S 5. $K$-Theory and the Trace Invariant \endsubhead
In this section we shall calculate the $K$-groups of the $C^{*}$-algebra 
$A:=\roman A_{\theta}^{5,3}(\alpha ,\beta ,\gamma ,\delta ,\epsilon 
)$
by means of the Pimsner-Voiculescu six term exact sequence [\PV ].  Since one of the groups
in the sequence turns out to have torsion elements, the application of this result 
requires careful examination.

\proclaim{6. Theorem} For the C*-algebra 
$\roman A_{\theta}^{5,3}(\alpha ,\beta ,\gamma ,\delta ,\epsilon )$,
one has $K_0=K_1=\Bbb Z^6\oplus \Bbb Z_{\alpha}$.
\endproclaim

\flushpar PROOF. To prove this theorem, we combine two applications 
of the PV sequence corresponding to two
presentations P1 and P2 of $A$ as follows.

\subhead P1 \endsubhead
In view of (CR), let $B_1=C^{*}(X,V,U)$ and let $\roman A\roman d_W$,
with
$$\roman A\roman d_W(X)=X,\qquad \roman A\roman
d_W(V)=\lambda^{-\delta}V,\qquad \roman A\roman d_W(U)=\lambda^{
-\epsilon}U,$$
generate an action of $\Bbb Z$ on $B_1$, so that $A=C^{*}(B_1,\Bbb Z
)$.
Applying the PV sequence to $B_1$, viewed as the crossed product of 
$C(\Bbb T^2)=C^{*}(X,V)$ by the automorphism $\roman A\roman d_U$, it is not hard to see that
$K_0(B_1)=\Bbb Z^3$ and $K_1(B_1)=\Bbb Z^3\oplus \Bbb Z_{\alpha}$.  Since $
\roman A\roman d_W$
is homotopic to the identity, the PV sequence immediately gives
$$K_1(A)=\Bbb Z^6\oplus \Bbb Z_{\alpha}.$$
However, since in the short exact sequence
$$\CD
0@>>>K_0(B_1)@>{i_{*}}>>K_0(A)@>{\delta}>>K_1(B_1)@>>>0\endCD
$$
$K_1(B_1)$ has torsion, we cannot readily obtain $K_0(A)$.  For this, the next
presentation will help.

\bigpagebreak

\subhead P2 \endsubhead
In view of (CR), we can also let $B_2=C^{*}(X,V,W)=C(\Bbb T)\otimes 
A_{\delta\theta}$,
where $C(\Bbb T)=C^{*}(X)$ and $A_{\delta\theta}=C^{*}(V,W)$.  Let $
\sigma =\roman A\roman d_U$, with
$$\sigma (X)=\lambda^{\gamma}X,\qquad\sigma (V)=\lambda^{\beta}X^{
\alpha}V,\qquad\sigma (W)=\lambda^{\epsilon}W,$$
generate an action of $\Bbb Z$ on $B_2$, so that $A=C^{*}(B_2,\Bbb Z
)$.
In this case the PV sequence becomes
$$\CD
K_0(B_2)@>{id_{*}-\sigma_{*}}>>K_0(B_2)@>{i_{*}}>>K_0(A)\\
@A{\delta_1}AA@.@VV{\delta_0}V\\
K_1(A)@<{i_{*}}<<K_1(B_2)@<{id_{*}-\sigma_{*}}<<K_1(B_2)\endCD
\tag *$$
It is not hard to see that a basis for $K_1(B_2)=\Bbb Z^4$ is given by 
$\{[X],[V],[W],[\xi ]\}$ where 
$\xi =X\otimes e+1\otimes (1-e)$ and $e=e(V,W)$ is a Rieffel projection in
$A_{\delta\theta}$ of trace $\delta\theta\mod 1$.  Also, a basis of $
K_0(B_2)=\Bbb Z^4$
is given by $\{[1],[e],B_{XV},B_{XW}\}$ where $B_{XV}=[P_{XV}]-[1]$ is
the Bott element in $X,V$ and $P_{XV}$ the usual Bott projection in the commuting
variables $X,V$.  The action of $id_{*}-\sigma_{*}$ on $K_1(B_2)$ is given by
$$id_{*}-\sigma_{*}:\qquad [X]\mapsto 0,\qquad [V]\mapsto -\alpha 
[X],\qquad [W]\mapsto 0,\qquad [\xi ]\mapsto m\alpha [X]$$
for some integer $m$, as shown by the following lemma.  The action of $
id_{*}-\sigma_{*}$
on $K_0(B_2)$ is given by
$$id_{*}-\sigma_{*}:\qquad [1]\mapsto 0,\qquad [e]\mapsto\alpha B_{
XW},\qquad B_{XW}\mapsto 0,\qquad B_{XV}\mapsto 0.$$
The action on $[e]$ is also shown in the following

\bigpagebreak

\proclaim{7. Lemma} We have $\sigma_{*}[e]=[e]-\alpha B_{XW}$ in $
K_0(B_2)$, and 
$\sigma_{*}[\xi ]=[\xi ]+m\alpha [X]$ for some integer $m$.
\endproclaim

\flushpar PROOF.
The proof of the first equality can be established using an argument quite similar 
to that of the proof of Lemma 4.2 of [\SW ].  Hence the kernel of $
id_{*}-\sigma_{*}$ 
on $K_0(B_2)$ is $\Bbb Z^3$.  For the second equality, let 
$\eta =(id_{*}-\sigma_{*})[\xi ]$.  From P1 and (*) we have
$$\Bbb Z^6\oplus \Bbb Z_{\alpha}=K_1(A)=\Bbb Z^3\oplus Im(i_{*})=\Bbb Z^
3\oplus\frac {K_1(B_2)}{Im(id_{*}-\sigma_{*})}=\Bbb Z^3\oplus\frac {
K_1(B_2)}{\Bbb Z\alpha [X]+\Bbb Z\eta}.$$
Thus
$$\frac {K_1(B_2)}{\Bbb Z\alpha [X]+\Bbb Z\eta}=\Bbb Z^3\oplus \Bbb Z_{\alpha}.\tag **$$
But since $K_1(B_2)=\Bbb Z^4$, it follows that the subgroup 
$\Bbb Z\alpha [X]+\Bbb Z\eta$ must have rank one.
\footnote{If $0\to F_1\to G\to F_2\oplus H\to 0$ is a short exact sequence of finitely
generated Abelian groups, where $F_1,F_2$ are free groups and $H$ is torsion, then 
$rank(G)=rank(F_1)+rank(F_2)$.  This can be seen from the
naturally obtained short exact sequence $0\to F_1\oplus F_2\to G\to 
H\to 0$, from which the result follows.  (If $G$ has rank greater than that of a subgroup
$K$, then $G/K$ contains a non-torsion element.)}
Therefore, $\Bbb Z\alpha [X]+\Bbb Z\eta =\Bbb Zd[X]$
for some integer $d$.  Substituting this into (**) one gets $d=\alpha$ and so
$\eta\in \Bbb Z\alpha [X]$.  \qed

\bigpagebreak

It now follows that in $K_1(B_2)$ one has 
$Im(id_{*}-\sigma_{*})=\Bbb Z\alpha [X]$ and that
$Ker(id_{*}-\sigma_{*})=\Bbb Z^3$ whether $m$ is zero or not.  Therefore, from
the exactness of (*) we obtain $Im(\delta_0)=\Bbb Z^3$ and hence by Lemma 7
$$K_0(A)=\Bbb Z^3\oplus Im(i_{*})=\Bbb Z^3\oplus\frac {K_0(B_2)}{I
m(id_{*}-\sigma_{*})}=\Bbb Z^3\oplus\frac {K_0(B_2)}{\Bbb Z\alpha 
B_{XW}}=\Bbb Z^6\oplus \Bbb Z_{\alpha}$$
which completes the proof of Theorem 6. \qed

\bigpagebreak

\subhead The Trace Invariant\endsubhead

\proclaim{8. Theorem} The range of the unique trace on 
$K_0(\roman A_{\theta}^{5,3}(\alpha ,\beta ,\gamma ,\delta ,\epsilon))$ is 
$\Bbb Z+\Bbb Z\rho\theta + \Bbb Z\gamma\delta\theta^2$ where 
$\rho =\gcd \{\gamma ,\delta ,\epsilon\}$.
\endproclaim

\flushpar Note that this agrees with the the trace invariant 
$\Bbb Z+\Bbb Z\theta + \Bbb Z\theta^2$ of the algebra $\roman
A_\theta^{5,3}$ as done in
[\SW], section 2, in the case $(\alpha ,\beta ,\gamma ,\delta
,\epsilon)=(1,0,1,1,0)$.

\smallpagebreak

\flushpar PROOF.
First we make an appropriate change of variables for the unitary generators of
the algebra $A=\roman A_{\theta}^{5,3}(\alpha ,\beta ,\gamma ,\delta ,\epsilon)$.
Referring back to the defining relations (CR), pick integers $a,b,c,d$ such that 
$b\delta+a\epsilon=0,\ ad-bc=1$, and let
$$
U' = U^aV^b, \qquad V' = U^cV^d.  
$$
Then the commutation relations (CR), with $W$ remaining the same and $X$ suitably scaled,
become
$$
\split
U'V' &= X^{\alpha}V'U', \quad U'X=\lambda^{a\gamma}XU', \quad
V'W=\lambda^{d\delta+c\epsilon}WV', \\ 
U'W &= WU',  \qquad  V'X=\lambda^{c\gamma}XV', \quad WX=XW
\endsplit
$$
Let $B=C^*(X,U',V')$.  It is isomorphic to the crossed product of 
$C^*(X,U') = A_{a\gamma\theta}$ by $\Bbb Z$ and automorphism $\roman
A\roman d_{V'}$.
An easy application of Pimsner's trace formula shows that 
$$
\vartau_*K_0(B) = \Bbb Z + \Bbb Za\gamma\theta + \Bbb Z c\gamma\theta =
\Bbb Z + \Bbb Z \gamma\theta,
$$
since $(a,c)=1$.
Next, it is not hard to see that an application of the Pimsner-Voiculescu sequence to the
above crossed product presentation of $B$ gives the basis $\{[X],[V'],[U'],[\xi]\}$ for
$K_1(B)$, where $[X]$ has order $\alpha$, $\xi = 1-e + ew^*{V'}^*e$ is a unitary in $B$,
$e$ is a Rieffel projection in $A_{a\gamma\theta}$ of trace
$(a\gamma\theta)\text{mod}\,1$, and
$w$ is a unitary in $A_{a\gamma\theta}$ such that ${V'}^*eV' = wew^*$ (which exists by 
Rieffel's Cancellation Theorem [\MR]).  The underlying connecting homomorphism
$\partial:K_1(B) \to K_0(A_{a\gamma\theta})$ gives $\partial[\xi]=[e]$ and 
$\partial[V']=[1]$, the usual basis of $K_0(A_{a\gamma\theta})$.

To apply Pimsner's trace formula, one calculates the usual ``determinant" on the 
aforementioned basis, since the kernel of $id_*-(\roman A\roman d_W)_*$ is all of $K_1(B)$ (since
$\roman A\roman d_W$ is homotopic to the identity).
It is easy to see that this determinant (whose values are in $\Bbb R/\vartau_*K_0(B)$)
on the elements $[X],[V'],[U']$ gives the respective values $1,(d\delta+c\epsilon)\theta,1$.
For the $\xi$, since now $\roman A\roman d_W$ fixes $A_{a\gamma\theta}$ (and in particular $e$ and $w$),
one obtains
$$
\roman A\roman d_W(\xi)\xi^* = (1-e + \lambda^{d\delta+c\epsilon} ew^*{V'}^*e) (1-e + eV'we)
= 1-e + \lambda^{d\delta+c\epsilon} e.
$$
Now a simple homotopy path connecting this element to 1 is just 
$t \mapsto 1-e + e^{2\pi i\theta (d\delta+c\epsilon)t} e$, and the corresponding 
determinant gives the value $(d\delta+c\epsilon)\theta \vartau(e)$.  Since 
$\vartau(e)=a\gamma\theta \mod 1$, the range of the trace is 
$$
\vartau_*K_0(A) = \Bbb Z + \Bbb Z\gamma\theta + \Bbb Z (d\delta+c\epsilon)\theta 
+ \Bbb Z \gamma a(d\delta+c\epsilon)\theta^2.
$$
Now $a(d\delta+c\epsilon) = ad\delta + ac\epsilon - c(b\delta+a\epsilon) = \delta$, and
similarly $-b(d\delta+c\epsilon) = \epsilon$, thus showing that 
$d\delta+c\epsilon = \gcd\{\delta,\epsilon\}$.  Therefore, one gets
$\vartau_*K_0(A) = \Bbb Z + \Bbb Z\gcd\{\gamma,\delta,\epsilon\}\theta +
\Bbb Z \gamma \delta \theta^2$.  \qed

\bigpagebreak

\subhead Discussion of Classification \endsubhead

Next, let us consider briefly the classification of the algebras
$\roman A_{\theta}^{5,3}(\alpha ,\beta ,\gamma ,\delta ,\epsilon )$.
First, it is easy to show that 
$\roman A_{\theta}^{5,3}(\alpha ,\beta ,\gamma ,\delta ,\epsilon)\cong 
\roman A_{-\theta}^{5,3}(\alpha ,\beta ,\gamma ,\delta ,\epsilon)$.
Second, we note that the simple quotients
$\roman A_{\theta}^{5,3}=\roman A_{\theta}^{5,3}(1,0,1,1,0)$ have been almost
completely classified in [\SW ]; specifically, they have been classified for all
non-quartic irrationals (which are those that are not zeros of any polynomial
of degree at most 4 with integer coefficients). But generally,
with $\lambda =e^{2\pi i\theta}$ for an irrational $\theta$, the operator equations
$$UV=\lambda^{\beta}X^{\alpha}VU,\,\,\,\,\,\,UX=\lambda^{\gamma}XU
,\,\,\,\,\,\,VW=\lambda^{\delta}WV\,\,\,\text{ and }\,\,\,\,UW=\lambda^{
\epsilon}WU,\leqno (\roman C\roman R)$$
for $\roman A_{\theta}^{5,3}(\alpha ,\beta ,\gamma ,\delta ,\epsilon
)$
can be modified by changing some of the variables, i.e., by substituting
$X_0=e^{2\pi i\theta\beta /\alpha}X$ and putting $\lambda_0=\lambda^{
\rho}$, where $\rho =\roman g\roman c\roman d\{\gamma ,\delta ,\epsilon
\}$, and then
$\gamma_0=\gamma /\rho$, $\delta_0=\delta /\rho$ and $\epsilon_0=\epsilon
/\rho$ with $\roman g\roman c\roman d\{\gamma_0,\delta_0,\epsilon_0\}=1$.
The equations (CR) become
$$\left\{\aligned
&UV=X_0^{\alpha}VU,\,\,\,\,\,\,UX_0=\lambda_0^{\gamma_0}X_0U,\,\,\,\,\,\,
VW=\lambda_0^{\delta_0}WV\,\,\,\text{ and }\\
&UW=\lambda_0^{\epsilon_0}WU\,\,\,\text{ with }\,\,\,\,
\roman g\roman c\roman d\{\gamma_0,\delta_
0,\epsilon_0\}=1,\endaligned
\right.\leqno (\roman C\roman R_0)
$$ 
which are the equations for 
$\roman A_{\rho\theta}^{5,3}(\alpha,0,\gamma_0,\delta_0,\epsilon_0)$, so
$$
\roman A_{\theta}^{5,3}(\alpha ,\beta ,\gamma ,\delta ,\epsilon) \cong 
\roman A_{\rho\theta}^{5,3}(\alpha ,0,\gamma_0,\delta_0,\epsilon_0)
$$
where $\roman g\roman c\roman d\{\gamma_0,\delta_0,\epsilon_0\}=1$.
This reduces the classification to the class of
algebras $\roman A_{\theta}^{5,3}(\alpha ,0,\gamma ,\delta ,\epsilon )$ where
$\roman g\roman c\roman d\{\gamma ,\delta ,\epsilon \}=1$.

If two such $C^*$-algebras $A_j = 
\roman A_{\theta_j}^{5,3}(\alpha_j ,0,\gamma_j ,\delta_j ,\epsilon_j)$, $j=1,2$, are 
isomorphic, where now $\rho_j=\roman g\roman c\roman d\{
\gamma_j ,\delta_j ,\epsilon_j\}=1$, what constraints must
hold between their respective parameters?  As we observed in Theorem 6, one must have 
$\alpha_1=\alpha_2$.  By Theorem 8, one has
$$
\Bbb Z+\Bbb Z\theta_1 + \Bbb Z\gamma_1\delta_1\theta_1^2 =
\Bbb Z+\Bbb Z\theta_2 + \Bbb Z\gamma_2\delta_2\theta_2^2.
$$ 
One can show that if one assumes that $\theta_j$ are non-quadratic irrationals, then 
these trace invariants are equal if, and only if, there is a matrix
$S\in\roman{GL}(2,\Bbb Z)$ such that
$$
\left(\matrix \theta_2 \\ \gamma_2\delta_2\theta_2^2 \endmatrix\right)
=
S \left(\matrix \theta_1 \\ \gamma_1\delta_1\theta_1^2 \endmatrix\right)
\text{mod}\, \left(\matrix \Bbb Z \\ \Bbb Z \endmatrix\right).
$$
Further, one can more easily show that if $\theta_j$ are non-quartic irrationals (i.e.,
not roots of polynomials over $\Bbb Z$ of degree at most four), then the trace
invariants are equal if, and only if, 
$$
\theta_2 = (\pm \theta_1)\text{mod}\, 1, \qquad \text{and} \qquad
\gamma_2\delta_2\theta_2^2 = (\pm \gamma_1\delta_1\theta_1^2 + m
\theta_1) \text{mod}\, 1,
$$
for some integer $m$.  An interesting special situation can be considered.  For example,
if one fixes $\theta$ (assumed non-quartic for simplicity) and varies the other parameters, 
then the above shows that $\gamma_1\delta_1 = \gamma_2\delta_2$ will follow from 
$A_1 \cong A_2$. At this point it is not clear if the parameters can be determined more precisely than
this. For example, is it possible, if $\theta$ is held fixed, that the equalities $\gamma_1=\gamma_2$ 
and $\delta_1=\delta_2$ could fail to hold? 
This is unclear.  However, the following heuristic argument 
(based only on canonical considerations) suggests that perhaps $\gamma_j ,\delta_j$ 
are uniquely determined.

\smallpagebreak

Let us attempt to apply a canonical transformation of the unitary generators of the form
$$
U_1=U^{a_1}V^{b_1}W^{c_1},\qquad V_1=U^{a_2}V^{b_2}W^{c_2},\qquad 
W_1=U^{a_3}V^{b_3}W^{c_3},
$$
in the hope of changing $\gamma ,\delta$, by working out the commutation
relations and ensuring that they are preserved.
(We have kept $X$ the same since it is the only auxiliary unitary
that can occur if one looks at the most general transformation of the form 
$U_1=U^{a_1}V^{b_1}W^{c_1}X^{d_1},~V_1=U^{a_2}V^{b_2}W^{c_2}X^{d_2}$ --- in fact,
the commutator $[U_1,V_1]$ is a scalar multiple of $X^{\alpha (a_1b_2-a_2b_1)}$.)
The 3 by 3 matrix $T$ with rows $a_j,b_j,c_j$ should have determinant $\pm 1$ for the
new unitaries to generate the same $C^*$-algebra.
The first relation in (CR) demands 
that $a_1b_2-a_2b_1=1$ so as to keep $X^{\alpha}$ the same.  Also, since $
VX=XV$ and $WX=XW$ must be preserved, we should have $V_1X=XV_1$ and $W_1X=XW_1$.  
However, it is easy to see that these imply that $a_2=0$ and $a_3=0$, respectively 
(since $U$ does not commute with $X$).  
Since the relation between $V_1$ and $W_1$ does not contain $X$, one must have
$a_2b_3-a_3b_2=0$, which is already satisfied.  Similarly, for the relation between
$U_1$ and $W_1$ one must have $a_1b_3-a_3b_1=0$.  But since $a_3=0$ this gives $b_3=0$.  
This means that the matrix $T$ is upper triangular with $1$ or $-1$ on its diagonal, 
hence the transformation is reduced to
$$
U_1=U^{\pm 1}V^{b_1}W^{c_1},\qquad V_1=V^{\pm 1}W^{c_2},\qquad 
W_1=W^{\pm 1}.
$$
(where the third $\pm$ here is independent of the first two, which should both be 
$1$ or both $-1$).  In view of this transformation, however, the new commutation relations
are now forced to take the following form
$$
U_1V_1=X^{\alpha}V_1U_1,\quad U_1X=\lambda^{\pm\gamma}XU_1,\quad 
V_1W_1=\lambda^{\pm\delta}W_1V_1,\quad U_1W_1=\lambda^{\delta b_1\pm \epsilon}W_1U_1,
$$
(and of course $V_1X=XV_1$, $W_1X=XW_1$, and after one rescales $X$).  These are exactly
in the same form as the relations (CR).  In particular, the integer parameters $\gamma$ and
$\delta$, since they are assumed to be positive, have remained unchanged.
(Also unchanged is
$\epsilon$, since it is, by $(*)$ of Theorem 1, smaller than $\delta$.)
This seems to suggest that $\gamma$ and $\delta$ (and hence also $\epsilon$) are 
uniquely determined in an isomorphism classification theorem.  The broad scope of this 
classification problem, however, must be left to another time; the fact that $\gamma$ 
and $\delta$ are not clearly singled out in the invariants considered here, but appear
mixed, seems to present an obstacle to the classification of these $C^{*}$-algebras.  
(The authors doubt that the Ext invariant of Brown-Douglas-Fillmore contains any more 
information, though they have not checked this in detail.)

\newpage


\enddocument